\documentclass[12pt]{article}
\usepackage{amsmath,amssymb,hyperref}
\newtheorem{theorem}{Theorem}
\newtheorem{proposition}{Proposition}
\begin{document}
\begin{center}
\large
\textbf{Pointwise Properties of Fourier-Jacobi-Forms}
\end{center}
\begin{center}
\textbf{Bert Koehler}
\end{center}
\begin{abstract}
\noindent
Jacobi-Forms can be decomposed as a linear combination of
Thetafunctions with modular forms as coefficients. It is shown that the space of these
coefficient modular forms of Fourier-Jacobi-Forms, which come from Siegel
cusp forms, has full rank in every point of the Satake boundary, if
the index is 1, the weight is sufficiently large and the Satake boundary point has
trivial stabilizer in $\Gamma_{n-1}$. This yields a local automorphic embedding of
the Siegel modular variety. Klingen-Poincare series are the main tool. Despite of
this richness it is proved that there are more Jacobi index 1 cusp forms than
Fourier-Jacobi index 1 cusp forms for all sufficiently large weights extending a result
of Dulinski.
\end{abstract}
\normalsize
\textbf{Introduction:} We start by fixing some notations: Let $\Gamma_n=Sp(n,\mathbb{Z})$
be the full symplectic group acting on the Siegel upper half space $\mathbb{H}_n$. 
For $Z\in\mathbb{H}_n$ let
\begin{eqnarray*}
Z=\left(
\begin{array}{cc}
\widehat{Z}&\widehat{z}\\
\widehat{z}&z_{nn}
\end{array}
\right)
\end{eqnarray*}
where $\widehat{Z}\in\mathbb{H}_{n-1}$ and $z_{nn}\in\mathbb{H}_1$. Then every
Siegel modular form $F\in[\Gamma_n,r]$ has a Fourier-Jacobi-decomposition \cite{1}
\begin{eqnarray*}
F(Z)=\Phi_0(F)(\widehat{Z})+\sum_{m=1}^{\infty}\Phi_m(F)(\widehat{Z},\widehat{z})e^{2\pi imz_{nn}}
\end{eqnarray*}
Here the coefficients $\Phi_m(F)(\widehat{Z},\widehat{z})$ are examples of
Jacobi-forms and at least since \cite{2} there is a lot of work 
on the inverse problem of lifting
Jacobi-forms to Siegel modular forms \cite{3},\cite{4}. In this article we want to show that
there are sufficiently many Fourier-Jacobi cusp forms in the sense that they provide a local
embedding of the Siegel modular variety but on the other hand are still not enough
to generate all Jacobi index 1 cusp forms extending the non-generation result of \cite{3}
to scalar index 1 case.\\
Every
Fourier-Jacobi-Form $\Phi_m(F)(\widehat{Z},\widehat{z})$ can be decomposed further
\begin{eqnarray*}
\Phi_m(F)(\widehat{Z},\widehat{z})=\sum_{a\in(\frac{1}{2m}\mathbb{Z}/\mathbb{Z})^{n-1}}
\phi_{m,a}(F)(\widehat{Z})\cdot\Theta_{m,a}(\widehat{Z},\widehat{z})
\end{eqnarray*}
Here the Theta-functions
\begin{eqnarray*}
\Theta_{m,a}(\widehat{Z},\widehat{z})=\sum_{k\in\mathbb{Z}^{n-1}}
e^{2\pi im(k+a)^t\widehat{Z}(k+a)+4\pi im(k+a)\widehat{z}}
\end{eqnarray*}
are a canonical basis of the space of sections of the corresponding line bundle
on the Abelian variety 
$mL\longrightarrow\mathbb{C}^{n-1}/\Lambda$ with
$\Lambda=Span_{\mathbb{Z}}(\widehat{E},\widehat{Z})$.\\
\\
For a point $\widehat{Z}\in\mathbb{H}_{n-1}$ let $Stab(\widehat{Z})\subset\Gamma_{n-1}$
be the finite set of transformations $\widehat{M}$ which fix $\widehat{Z}=\widehat{M}(\widehat{Z})$.
The set of points $\widehat{Z}$ with $Stab(\widehat{Z})\neq\{id\}$ is a high codimensional analytic
subset of $\mathbb{H}_{n-1}$. Let
\begin{eqnarray*}
\widetilde{\mathbb{H}}_{n-1}=\{\widehat{Z}\in\mathbb{H}_{n-1}\mbox{ with }
Stab(\widehat{Z})=\{id\}\}
\end{eqnarray*}
$\mathbb{H}_{n-1}/\Gamma_{n-1}$ is a quasiprojective algebraic variety \cite{5} and
$\widetilde{\mathbb{H}}_{n-1}/\Gamma_{n-1}\subset\mathbb{H}_{n-1}/\Gamma_{n-1}$ 
is an open, dense quasiprojective algebraic subvariety. For a weight $r$ let
$F_{r,1},\cdots,F_{r,N(r)}\in[\Gamma_n,r]_0$ be a basis of Siegel cusp forms of weight $r$
and let $\phi_{1,a}(F_{r,j})$, $1\leq j\leq N(r)$, $a\in(\frac{1}{2}\mathbb{Z}/\mathbb{Z})^{n-1}$
be the modular forms in the Fourier-Jacobi-decomposition of index 1.
\begin{theorem}
There is a sufficiently large weight $r>>2n+1$ such that for all
$\widehat{Z}\in\widetilde{\mathbb{H}}_{n-1}$ we have
\begin{eqnarray*}
rank\left(\phi_{1,a}(F_{r,j})(\widehat{Z})\right)_{1\leq j\leq N(r),
a\in(\frac{1}{2}\mathbb{Z}/\mathbb{Z})^{n-1}}=\max=2^{n-1}
\end{eqnarray*}
\end{theorem}
\textbf{Proof:} The set of points $[\widehat{Z}]\in\widetilde{\mathbb{H}}_{n-1}/
\Gamma_{n-1}$ with maximal rank is a Zariski open (maybe empty) subset.
If we can find for every $\widehat{Z}\in\widetilde{\mathbb{H}}_{n-1}$ a weight 
$r=r(\widehat{Z})$
with $rank=2^{n-1}$, then due to $\widetilde{\mathbb{H}}_{n-1}/\Gamma_{n-1}$
quasiprojective every chain of subvarieties is finite and so there is a common weight
$r$, which works for all $\widehat{Z}\in\widetilde{\mathbb{H}}_{n-1}$.\\
Let $\mathcal{F}_{n-1}\subset\mathbb{H}_{n-1}$ be the Siegel fundamental domain
and we can assume in the following that $\widehat{Z}\in\mathcal{F}_{n-1}\cap
\widetilde{\mathbb{H}}_{n-1}$ which implies for every $\widehat{M}\in\Gamma_{n-1}$
\begin{eqnarray*}
|\det(\widehat{C}\widehat{Z}+\widehat{D})|\geq 1\mbox{ for }
\widehat{M}=\left(
\begin{array}{cc}
\widehat{A}&\widehat{B}\\
\widehat{C}&\widehat{D}
\end{array}
\right)
\end{eqnarray*}
In the following let $T$ be a symmetric, positive definite $n\times n$-matrix which
is half-even, this means $T_{ii}\in\mathbb{N}$ and $T_{ij}\in\frac{1}{2}\mathbb{Z}$.
According to \cite{6} for every such matrix $T$ and for $r>2n+1$ 
the Klingen-Poincare-series
\begin{eqnarray*}
F(T,Z)=\sum_{[M]\in\Gamma_n/\{\mbox{translations}\}}(\det(CZ+D))^{-r}e^{2\pi i\sigma(TM(Z))}
\end{eqnarray*}
converges absolutely and uniformly in $\mathcal{F}_n$ and 
is a Siegel cusp form in $[\Gamma_n,r]_0$. By abuse of notation let
\begin{eqnarray*}
\Phi_1(F_T)(\widehat{Z},\widehat{z})=
\sum_{a\in(\frac{1}{2}\mathbb{Z}/\mathbb{Z})^{n-1}}
\phi(a,T)(\widehat{Z})\cdot\Theta_a(\widehat{Z},\widehat{z})
\end{eqnarray*}
be the Fourier-Jacobi decomposition of index 1 of $F(T,Z)$.
We want to find at least $2^{n-1}$ matrices $T=T(t)$ such that the matrix of
coefficient modular forms $(\phi(a,T(t))(\widehat{Z}))_{a,t}$ has full rank $2^{n-1}$.\\ 
Now $F(T,Z)$ can be written as
\begin{eqnarray*}
F(T,Z)=\sum_{[M]\in\Gamma_n/\Gamma_{n,0}}(\det(CZ+D))^{-r}\sum_{U\in GL(n,\mathbb{Z})}
e^{2\pi i\sigma(TU^t M(Z)U)}\\
=G(T,Z)+\sum_{[id]\neq[M]\in\Gamma_n/\Gamma_{n,0}}(\det(CZ+D))^{-r}
\sum_{U\in GL(n,\mathbb{Z})}e^{2\pi i\sigma(TU^t M(Z)U)}
\end{eqnarray*}
where
\begin{eqnarray*}
G(T,Z)=\sum_{U\in GL(n,\mathbb{Z})}e^{2\pi i\sigma(TU^t ZU)}=
\sum_{U\in GL(n,\mathbb{Z})}e^{2\pi i\sigma(UTU^t Z)}=\\
=\sum_{k=1}^{\infty}g_k(T,\widehat{Z},\widehat{z})e^{2k\pi iz_{nn}}
\end{eqnarray*}
We calculate the trace
\begin{eqnarray*}
\sigma(TU^t ZU)=z_{nn}\sum_{i,j=1}^n T_{ij}U_{nj}U_{ni}+
2\sum_{k=1}^{n-1}\widehat{z}_k\sum_{i,j=1}^n T_{ij}U_{kj}U_{ni}+
\sum_{k,l=1}^{n-1}\widehat{Z}_{kl}\sum_{i,j=1}^n T_{ij}U_{kj}U_{li}
\end{eqnarray*}
In the following we always choose
\begin{eqnarray*}
T=\left(
\begin{array}{cc}
\widehat{T}&T_{in}\\
T_{in}^t&T_{nn}
\end{array}
\right)
\end{eqnarray*}
with $T_{nn}=1$, $T_{in}=0$ for $2\leq i\leq n-1$, $T_{1n}\in\{0,\frac{1}{2}\}$
and $\widehat{T}$ with eigenvalues $\geq 3$. 
For such $T$ we have for the coefficient of $z_{nn}$ in the above trace
\begin{eqnarray*}
\sum_{i,j=1}^n T_{ij}U_{nj}U_{ni}=\sum_{i,j=1}^{n-1}\widehat{T}_{ij}U_{ni}U_{nj}+
2T_{1n}U_{n1}U_{nn}+U_{nn}^2\geq\\
3\sum_{i=1}^{n-1}U_{ni}^2-|U_{n1}||U_{nn}|+U_{nn}^2\geq
2\sum_{i=1}^{n-1}U_{ni}^2+\frac{1}{2}U_{nn}^2
\end{eqnarray*}
So the coefficient of $z_{nn}$ becomes 1
precisely if $U_{ni}=0$ for all $1\leq i\leq n-1$ and $U_{nn}=\pm 1$. In this case unimodular
$U$ is of the form
\begin{eqnarray*}
U=\left(
\begin{array}{cc}
\widehat{U}&U_{in}\\
0&\pm 1
\end{array}
\right)
\end{eqnarray*}
with $\widehat{U}\in GL(n-1,\mathbb{Z})$ unimodular and $U_{in}\in\mathbb{Z}$ with 
$1\leq i\leq n-1$ can 
vary free. The coefficients of $\widehat{z}_k$ then reduce to
\begin{eqnarray*}
2\sum_{i,j=1}^{n}T_{ij}U_{kj}U_{ni}=\pm 2\sum_{j=1}^n T_{nj}U_{kj}=
\pm 2(U_{kn}+T_{1n}U_{k1})
\end{eqnarray*}
So we see
\begin{eqnarray*}
g_1(T,\widehat{Z},\widehat{z})=\sum_{\widehat{U}\in GL(n-1,\mathbb{Z})}\sum_{
(U_{1n},...,U_{n-1,n})\in\mathbb{Z}^{n-1}}e^{2\pi i(\pm 2\sum_{k=1}^{n-1}
\widehat{z}_k(U_{kn}+T_{1n}U_{k1})}\cdot\\
\exp\left(2\pi i\left(\sum_{k,l=1}^{n-1}\widehat{Z}_{kl}
(\sum_{i,j=1}^{n-1}T_{ij}U_{kj}U_{li}+T_{1n}(U_{kn}U_{l1}+U_{k1}U_{ln})+
U_{kn}U_{ln})\right)\right)
\end{eqnarray*}
In the following let $a\in\left(\frac{1}{2}\mathbb{Z}/\mathbb{Z}\right)^{n-1}$ 
be a theta-characteristic for $m=1$ and
\begin{eqnarray*}
\Theta_{a}(\widehat{Z},\widehat{z})=\sum_{g\in\mathbb{Z}^{n-1}}
e^{2\pi i(g+a)^t\widehat{Z}(g+a)+4\pi i(g+a)^t\widehat{z}}
\end{eqnarray*}
the corresponding theta-function. Let
\begin{eqnarray*}
GL_a(n-1,\mathbb{Z})=\{\widehat{U}\in GL(n-1,\mathbb{Z})\mbox{ with }
(\widehat{U})_1=
(\widehat{U}_{11},...,\widehat{U}_{n-1,1})\equiv 2a\mbox{ mod }2\}
\end{eqnarray*}
Then we can write
\begin{eqnarray*}
g_1(T,\widehat{Z},\widehat{z})=\sum_{\widehat{U}\in GL(n-1,\mathbb{Z})}
e^{2\pi i\sigma(\widehat{U}\widehat{T}\widehat{U}^t\widehat{Z})-
2\pi iT_{1n}^2(\widehat{U})_1^t\widehat{Z}(\widehat{U})_1}\cdot\\
\sum_{g\in\mathbb{Z}^{n-1}}e^{2\pi i(g+T_{1n}(\widehat{U})_1)^t\widehat{Z}
(g+T_{1n}(\widehat{U})_1)\pm 4\pi i(g+T_{1n}(\widehat{U})_1)^t\widehat{z}}=
\end{eqnarray*}
\begin{eqnarray*}
=\sum_{a\in\left(\frac{1}{2}\mathbb{Z}/\mathbb{Z}\right)^{n-1}}
\sum_{\widehat{U}\in GL_a(n-1,\mathbb{Z})}
e^{2\pi i\sigma(\widehat{U}\widehat{T}\widehat{U}^t\widehat{Z})-
2\pi iT_{1n}^2(\widehat{U})_1^t\widehat{Z}(\widehat{U})_1}\cdot\\
\sum_{g\in\mathbb{Z}^{n-1}}e^{2\pi i(g+2aT_{1n})^t\widehat{Z}
(g+2aT_{1n})\pm 4\pi i(g+2aT_{1n})^t\widehat{z}}
\end{eqnarray*}
So if we choose $T_{1n}=\frac{1}{2}$ then we get
\begin{eqnarray*}
g_1(T,\widehat{Z},\widehat{z})=
\sum_{a\in\left(\frac{1}{2}\mathbb{Z}/\mathbb{Z}\right)^{n-1},a\neq(0,...,0)}
\psi(a,T)(\widehat{Z})\cdot\Theta_{a}(\widehat{Z},\widehat{z})
\end{eqnarray*}
with
\begin{eqnarray*}
\psi(a,T)(\widehat{Z})=\sum_{\widehat{U}\in GL_a(n-1,\mathbb{Z})}
e^{2\pi i\sigma(\widehat{U}\widehat{T}\widehat{U}^t\widehat{Z})-
\frac{1}{2}\pi i(\widehat{U})_1^t\widehat{Z}(\widehat{U})_1}
\end{eqnarray*}
and if we choose $T_{1n}=0$ then
\begin{eqnarray*}
g_1(T,\widehat{Z},\widehat{z})=\psi((0,0,...,0),T)(\widehat{Z})\cdot
\Theta_{(0,...,0)}(\widehat{Z},\widehat{z})
\end{eqnarray*}
with
\begin{eqnarray*}
\psi((0,0,...,0),T)(\widehat{Z})=\sum_{\widehat{U}\in GL(n-1,\mathbb{Z})}
e^{2\pi i\sigma(\widehat{U}\widehat{T}\widehat{U}^t\widehat{Z})}
\end{eqnarray*}
If we can find 
$t=2,...,2^{n-1}$ matrices $T(t)$
with $T(t)_{1n}=\frac{1}{2}$ such that the 
$(2^{n-1}-1)\times(2^{n-1}-1)$-matrix
$(\psi(a,T(t))(\widehat{Z}))_{a\neq(0,0,...,0),t}$ has rank $2^{n-1}-1$ and if
we can find a further $T(1)$ with $T(1)_{1n}=0$ and
$\psi(a=(0,0,...,0),T(1))(\widehat{Z})\neq 0$, then the Fourier-Jacobi decomposition
of $g_1(T(t),\widehat{Z},\widehat{z})$, $1\leq t\leq 2^{n-1}$ has full rank $2^{n-1}$ in
$\widehat{Z}$.\\
\\
Let $\widehat{Z}=\widehat{X}+i\widehat{Y}$ and
\begin{eqnarray*}
0<\varepsilon_0=\min\{\sigma(\widehat{U}^t\widehat{Y}\widehat{U}),
\widehat{U}\in GL(n-1,\mathbb{Z})\}
\end{eqnarray*}
Let
\begin{eqnarray*}
H=\{\widehat{U}\in GL(n-1,\mathbb{Z})
\mbox{ with }\widehat{U}^t\widehat{Y}\widehat{U}=\widehat{Y}\}
\end{eqnarray*}
be the unimodular stabilizer of $\widehat{Y}$. This finite group can be nontrivial despite
the assumption that $\widehat{Z}$ has trivial stabilizer. Let
\begin{eqnarray*}
H=\{\widehat{V}_1,...,\widehat{V}_{2N}\}
\end{eqnarray*}
where $\widehat{V}_1=\widehat{E}$ and $\widehat{V}_2=-\widehat{E}$. If $N>1$ then
we have $\widehat{V}_j^t\widehat{X}\widehat{V}_j\neq\widehat{X}$ for $j>2$. Let
\begin{eqnarray*}
\widetilde{H}=\{\widehat{U}\in GL(n-1,\mathbb{Z})
\mbox{ with }\sigma(\widehat{U}^t\widehat{Y}\widehat{U})=\varepsilon_0\}
\end{eqnarray*}
which is a finite set. $H$ acts on $\widetilde{H}$ from left and so let
$\widehat{U}_1,...,\widehat{U}_M$ be nonequivalent representatives for $\widetilde{H}$, so
$\widehat{U}_i^t\widehat{Y}\widehat{U}_i\neq\widehat{U}_j^t\widehat{Y}\widehat{U}_j$ 
for $i\neq j$. 
Then we can find a rational symmetric matrix $\widetilde{T}$ sufficiently close
to the identity matrix $\widehat{E}$ with
\begin{eqnarray*}
\sigma(\widetilde{T}\widehat{U}_1^t\widehat{Y}\widehat{U}_1)<
\sigma(\widetilde{T}\widehat{U}_2^t\widehat{Y}\widehat{U}_2)<\cdots<
\sigma(\widetilde{T}\widehat{U}_M^t\widehat{Y}\widehat{U}_M)
\end{eqnarray*}
and
\begin{eqnarray*}
\sigma(\widetilde{T}\widehat{U}_M^t\widehat{Y}\widehat{U}_M)<
\sigma(\widetilde{T}\widehat{U}^t\widehat{Y}\widehat{U})\mbox{ for all }
\widehat{U}\in GL(n-1,\mathbb{Z})\backslash\widetilde{H}
\end{eqnarray*}
For a certain $q_0\in\mathbb{N}$ we set $\widehat{T}=q_0\widetilde{T}$ with 
$\widehat{T}$ symmetric, integer-valued and large (positive) eigenvalues.
Let
\begin{eqnarray*}
\varepsilon_1=\sigma(\widehat{T}\widehat{U}_1^t\widehat{Y}\widehat{U}_1)>0
\mbox{ and }
\varepsilon_2=\sigma(\widehat{T}\widehat{U}_2^t\widehat{Y}\widehat{U}_2)=
\varepsilon_1+\delta>\varepsilon_1
\end{eqnarray*}
then we have by construction
$\sigma(\widehat{T}\widehat{U}^t\widehat{Y}\widehat{U})\geq
\varepsilon_1+\delta$ for all $\widehat{U}\in GL(n-1,\mathbb{Z})$ with
$\widehat{U}\notin\{\widehat{V}_j\widehat{U}_1\mbox{ }j=1,...,2N\}$.
Because $\widehat{V}_j^t\widehat{X}\widehat{V}_j\neq\widehat{X}$ mod $\mathbb{Z}$
for $j>2$ we can assume
\begin{eqnarray*}
e^{2\pi i\sigma(\widehat{T}\widehat{U}_1^t\widehat{V}_i^t
\widehat{X}\widehat{V}_i\widehat{U}_1)}
\neq
e^{2\pi i\sigma(\widehat{T}\widehat{U}_1^t\widehat{V}_j^t
\widehat{X}\widehat{V}_j\widehat{U}_1)}\mbox{ unless }\widehat{V}_i=\pm\widehat{V}_j
\end{eqnarray*}
Now we consider $q\widehat{T}$ for $q\in\mathbb{N}$. 
For $q\longrightarrow\infty$ all terms in
\begin{eqnarray*}
\psi((0,0,...,0),qT)(\widehat{Z})=\sum_{\widehat{U}\in GL(n-1,\mathbb{Z})}
e^{2\pi iq\sigma(\widehat{T}\widehat{U}^t\widehat{Z}\widehat{U})}
\end{eqnarray*}
are dominated by
\begin{eqnarray*}
e^{2\pi iq\sigma(\widehat{T}\widehat{U}_1^t\widehat{Y}\widehat{U}_1)}\cdot
\sum_{j=1}^{2N}
e^{2\pi iq\sigma(\widehat{T}\widehat{U}_1^t\widehat{V}_j^t\widehat{X}
\widehat{V}_j\widehat{U}_1)}=
2e^{2\pi iq\sigma(\widehat{T}\widehat{U}_1^t\widehat{Y}\widehat{U}_1)}\cdot
\sum_{j=2l-1}
e^{2\pi iq\sigma(\widehat{T}\widehat{U}_1^t\widehat{V}_j^t\widehat{X}
\widehat{V}_j\widehat{U}_1)}
\end{eqnarray*}
Because of
$e^{2\pi i\sigma(\widehat{T}\widehat{U}_1^t\widehat{V}_i^t
\widehat{X}\widehat{V}_i\widehat{U}_1)}
\neq
e^{2\pi i\sigma(\widehat{T}\widehat{U}_1^t\widehat{V}_j^t
\widehat{X}\widehat{V}_j\widehat{U}_1)}$ the function
\begin{eqnarray*}
s\longmapsto\sum_{j=2l-1}
e^{2\pi i\sigma(\widehat{T}\widehat{U}_1^t\widehat{V}_j^t\widehat{X}
\widehat{V}_j\widehat{U}_1)}\cdot
\left(1-s\cdot e^{2\pi i\sigma(\widehat{T}\widehat{U}_1^t\widehat{V}_j^t\widehat{X}
\widehat{V}_j\widehat{U}_1)}\right)^{-1}
\end{eqnarray*}
has $N$ different poles on the unit circle and so radius of convergence is 1. So
there are infinitely many $q\in\mathbb{N}$ with
\begin{eqnarray*}
\left|\sum_{j=2l-1}
e^{2\pi iq\sigma(\widehat{T}\widehat{U}_1^t\widehat{V}_j^t\widehat{X}
\widehat{V}_j\widehat{U}_1)}\right|>e^{-\pi q\delta}
\end{eqnarray*}
This implies $\psi((0,0,...,0),qT)(\widehat{Z})\neq 0$ for infinitely many $q$
sufficiently large.\\
Now for $a\neq(0,...,0)$ we first fix the same $\widehat{T}$ combined with $T_{1n}=\frac{1}{2}$.
Let $\widehat{W}_1,...,\widehat{W}_{2^{n-1}-1}$, 
$\widehat{W}_1=\widehat{E}$ be unimodular matrices such that
$\frac{1}{2}(\widehat{U}_1\widehat{W}_j)_1$ runs through all theta-characteristics 
$\neq(0,...,0)$ for $j=1,...,2^{n-1}-1$. Let $\widehat{T}_j=
\widehat{W}_j^{-1}\widehat{T}(\widehat{W}_j^{-1})^t$ which is symmetric, 
integer valued and positive. Then by construction we have
\begin{eqnarray*}
\sigma(\widehat{T}_j\widehat{W}_j^t\widehat{U}_1^t\widehat{Y}
\widehat{U}_1\widehat{W}_j)=
\sigma(\widehat{T}\widehat{U}_1^t\widehat{Y}\widehat{U}_1)=\min(
\sigma(\widehat{T}_j\widehat{U}^t\widehat{Y}\widehat{U}),\mbox{ }
\widehat{U}\in GL(n-1,\mathbb{Z}))
\end{eqnarray*}
and with all these matrices $\widehat{T}=\widehat{T}_1,...,\widehat{T}_{2^{n-1}-1}$
the dominant terms are shifted to every theta characteristic $\neq(0,...,0)$.
As before let $q\in\mathbb{N}$ and consider the corresponding coefficient functions
\begin{eqnarray*}
\psi(a,q\widehat{T}_j,T_{1n})(\widehat{Z})=\sum_{\widehat{U}\in GL_a(n-1,\mathbb{Z})}
e^{2\pi i\sigma(q\widehat{T}_j\widehat{U}^t\widehat{Z}\widehat{U})-
\frac{1}{2}\pi i(\widehat{U})_1^t\widehat{Z}(\widehat{U})_1}
\end{eqnarray*}
Assume by contradiction that there are complex numbers $c_a=c_a(\widehat{Z})$ 
not all 0 such that for all $q$ and all $j=1,...,2^{n-1}-1$ we have
\begin{eqnarray*}
\sum_{a\neq(0,...,0)}c_a\psi(a,q\widehat{T}_j,T_{1n})(\widehat{Z})=0
\end{eqnarray*}
This is equivalent to
\begin{eqnarray*}
\sum_{a\neq(0,...,0)}c_a\left(\sum_{\widehat{U}\in GL_a(n-1,\mathbb{Z})}
e^{2\pi i\sigma(\widehat{T}_j\widehat{U}^t\widehat{Z}\widehat{U})-
\frac{1}{2}\pi i(\widehat{U})_1^t\widehat{Z}(\widehat{U})_1}\cdot
\left(1-s\cdot
e^{2\pi i\sigma(\widehat{T}_j\widehat{U}^t\widehat{Z}\widehat{U})}\right)^{-1}\right)=0
\end{eqnarray*}
for all $j=1,...,2^{n-1}-1$ and all complex numbers $s$. Now define theta-characteristics
$a(j,k)$ with $k=1,...,N$ by
\begin{eqnarray*}
2a(j,k)\equiv(\widehat{V}_{2k-1}\widehat{U}_1\widehat{W}_j)_1\mbox{ mod }2
\end{eqnarray*}
Then the partial sum of
 terms in the above series with poles closest to the unit circle must be 0 and this is given by
\begin{eqnarray*}
\sum_{k=1}^N c_{a(j,k)}\Big(e^{2\pi i\sigma(\widehat{T}_j\widehat{W}_j^t
\widehat{U}_1^t\widehat{V}_{2k-1}^t\widehat{Z}
\widehat{V}_{2k-1}\widehat{U}_1\widehat{W}_j)-
\frac{1}{2}\pi i(\widehat{V}_{2k-1}\widehat{U}_1\widehat{W}_j)_1^t
\widehat{Z}(\widehat{V}_{2k-1}\widehat{U}_1\widehat{W}_j)_1}\cdot\\
\Big(1-s\cdot
e^{2\pi i\sigma(\widehat{T}_j\widehat{W}_j^t
\widehat{U}_1^t\widehat{V}_{2k-1}^t\widehat{Z}
\widehat{V}_{2k-1}\widehat{U}_1\widehat{W}_j)}\Big)^{-1}\Big)=0
\end{eqnarray*}
Because of
\begin{eqnarray*}
e^{2\pi i\sigma(\widehat{T}_j\widehat{W}_j^t
\widehat{U}_1^t\widehat{V}_{2l-1}^t\widehat{X}
\widehat{V}_{2l-1}\widehat{U}_1\widehat{W}_j)}\neq
e^{2\pi i\sigma(\widehat{T}_j\widehat{W}_j^t
\widehat{U}_1^t\widehat{V}_{2k-1}^t\widehat{X}
\widehat{V}_{2k-1}\widehat{U}_1\widehat{W}_j)}\mbox{ for }l\neq k
\end{eqnarray*}
we conclude $c_{a(j,k)}=0$ for all $k$ and all $j$ and so all $c_a=0$. So there
are $2^{n-1}-1$ values of $q=q_j$ such that
\begin{eqnarray*}
\det\left(\psi(a,q_j\widehat{T}_j,T_{1n})(\widehat{Z})\right)_{a\neq(0,...,0),
j=1,...,2^{n-1}-1}\neq 0
\end{eqnarray*}
So we have found $2^{n-1}$ matrices $T_j$ such that the Theta-decomposition of
$(G(T_j,Z))_{1\leq j\leq 2^{n-1}}$ has full rank $2^{n-1}$ where we remind 
\begin{eqnarray*}
F(T,Z)=\sum_{[M]\in\Gamma_n/\Gamma_{n,0}}(\det(CZ+D))^{-r}\sum_{U\in GL(n,\mathbb{Z})}
e^{2\pi i\sigma(TU^t M(Z)U)}\\
=G(T,Z)+\sum_{[id]\neq[M]\in\Gamma_n/\Gamma_{n,0}}(\det(CZ+D))^{-r}
\sum_{U\in GL(n,\mathbb{Z})}e^{2\pi i\sigma(TU^t M(Z)U)}
\end{eqnarray*}
The second series must be decomposed further. Let
\begin{eqnarray*}
\Gamma_{n,n-1}=\{M=\left(
\begin{array}{cc}
A&B\\
C&D
\end{array}
\right)\in\Gamma_n
\mbox{ with }C_{in}=0\mbox{ for all }1\leq i\leq n\}
\end{eqnarray*}
Every transformation $M\in\Gamma_{n,n-1}$ is unimodular equivalent to a reduced form
where
\begin{eqnarray*}
A=\left(
\begin{array}{cc}
\widehat{A}&0\\
0&1
\end{array}
\right)\quad
C=\left(
\begin{array}{cc}
\widehat{C}&0\\
0&0
\end{array}
\right)\quad
D=\left(
\begin{array}{cc}
\widehat{D}&0\\
0&1
\end{array}
\right)
\end{eqnarray*}
and $\widehat{M}\in\Gamma_{n-1}$. Now for a generic point 
$\widehat{Z}\in\mathcal{F}_{n-1}$ we have 
$|\det(\widehat{C}\widehat{Z}+\widehat{D})|>1$
for all $\widehat{M}$ with $\widehat{C}\neq (0)$. For 
$M\in\Gamma_n\backslash\Gamma_{n,n-1}$, so with at least one index $i$ with
$C_{in}\neq 0$, we will show at the end
\begin{eqnarray*}
|\det(CZ+D)|\geq Const(C)\cdot y_{nn}>2\mbox{ for all }y_{nn}>Const
\mbox{ sufficiently large}
\end{eqnarray*}
This means that for a generic point $\widehat{Z}\in\mathcal{F}_{n-1}$
all terms in
\begin{eqnarray*}
\sum_{[id]\neq[M]\in\Gamma_n/\Gamma_{n,0}}(\det(CZ+D))^{-r}
\sum_{U\in GL(n,\mathbb{Z})}e^{2\pi i\sigma(TU^t M(Z)U)}
\end{eqnarray*}
can be made arbitrary small by choosing the weight $r$
very large. Now the modular forms in the Theta-decomposition can be retained by
\begin{eqnarray*}
\phi_{1,a}(F)(\widehat{Z})e^{-2\pi y_{nn}}\Theta_{1,a}(2\widehat{Z},0)=
\int_0^1\int_0^1...\int_0^1 F(Z)e^{-2\pi ix_{nn}}
\Theta_{1,a}(\widehat{Z},-\widehat{z})dx_{1n}...dx_{nn}
\end{eqnarray*}
So for sufficiently large $r$ all transformations $M\in\Gamma_n$ with
$C\neq(0)$ do not change the rank of the Theta-decomposition of
$(G(T_j,Z))_{1\leq j\leq2^{n-1}}$, if $\widehat{Z}$ is generic in the sense that
$|\det(\widehat{C}\widehat{Z}+\widehat{D})|>1$ for all $\widehat{C}\neq(0)$.
This proves the Theorem for generic $\widehat{Z}$.\\
\\
In case $\widehat{Z}\in\mathcal{F}_{n-1}$ is not generic, so
$|\det(\widehat{C}_p\widehat{Z}+\widehat{D}_p)|=1$ for some 
finitely many transformations $M_p\in\Gamma_{n,n-1}$, $p=2,...,P$
with $\widehat{C}_p\neq(0)$ we can extend the former argument (we let
$M_1=id$). We have
\begin{eqnarray*}
F(T,Z)=\sum_{[M]\in\Gamma_n/\Gamma_{n,0}}(\det(CZ+D))^{-r}\sum_{U\in GL(n,\mathbb{Z})}
e^{2\pi i\sigma(TU^t M(Z)U)}\\
=\sum_{p=1}^P
(\det(\widehat{C}_p\widehat{Z}+\widehat{D}_p))^{-r}
\sum_{U\in GL(n,\mathbb{Z})}e^{2\pi i\sigma(TU^t M_p(Z)U)}+\\
\sum_{[M]\neq[id],[M_1],...,[M_p]}
(\det(CZ+D))^{-r}\sum_{U\in GL(n,\mathbb{Z})}
e^{2\pi i\sigma(TU^t M(Z)U)}
\end{eqnarray*}
The contribution of the second series to the rank can be made small for $r$ large and so
can be neglected. Because $M_p\in\Gamma_{n,n-1}$ we have
\begin{eqnarray*}
M_p(Z)=\left(
\begin{array}{cc}
\widehat{M}_p(\widehat{Z})&(\widehat{Z}\widehat{C}_p^t+\widehat{D}_p^t)^{-1}\widehat{z}\\
\widehat{z}^t(\widehat{C}_p\widehat{Z}+\widehat{D}_p)^{-1}&
z_{nn}-\widehat{z}^t(\widehat{C}_p\widehat{Z}+\widehat{D}_p)^{-1}\widehat{C}_p
\widehat{z}
\end{array}
\right)
\end{eqnarray*}
So if we choose symmetric half-integer $T>0$ as before with $T_{nn}=1$ and
$\widehat{T}>3\widehat{E}$ then exactly those unimodular $U$ contribute to
the first Fourier-coefficient in $z_{nn}$ with
\begin{eqnarray*}
U=\left(
\begin{array}{cc}
\widehat{U}&U_{in}\\
0&\pm 1
\end{array}
\right)
\end{eqnarray*}
as before and so all terms of
\begin{eqnarray*}
\sum_{U\in GL(n,\mathbb{Z})}e^{2\pi i\sigma(TU^t M_p(Z)U)}
\end{eqnarray*}
contributing to the first Fourier-coefficient in $z_{nn}$ sum up to
\begin{eqnarray*}
g_1(T,\widehat{Z}_p,\widehat{z}_p)=
\sum_{a\in\left(\frac{1}{2}\mathbb{Z}/\mathbb{Z}\right)^{n-1},a\neq(0,...,0)}
\kappa_p\cdot\psi(a,T)(\widehat{Z}_p)\cdot\Theta_{a}(\widehat{Z}_p,\widehat{z}_p)
\end{eqnarray*}
for $T_{1n}=\frac{1}{2}$ and
\begin{eqnarray*}
g_1(T,\widehat{Z}_p,\widehat{z}_p)=\kappa_p\cdot\psi((0,0,...,0),T)(\widehat{Z}_p)\cdot
\Theta_{(0,...,0)}(\widehat{Z}_p,\widehat{z}_p)
\end{eqnarray*}
for $T_{1n}=0$. Here is $\widehat{Z}_p=\widehat{M}_p(\widehat{Z})$,
$\widehat{z}_p=(\widehat{Z}\widehat{C}_p^t+\widehat{D}_p^t)^{-1}\widehat{z}$ and
\begin{eqnarray*}
\kappa_p=e^{-2\pi i\widehat{z}^t
(\widehat{C}_p\widehat{Z}+\widehat{D}_p)^{-1}\widehat{C}_p
\widehat{z}}
\end{eqnarray*}
Now we invoke the Theta-transformation-formula
\begin{eqnarray*}
\Theta_a(\widehat{Z}_p,\widehat{z}_p)=
\sqrt{\det(\widehat{C}_p\widehat{Z}+\widehat{D}_p)}\cdot
e^{2\pi i\widehat{z}^t
(\widehat{C}_p\widehat{Z}+\widehat{D}_p)^{-1}\widehat{C}_p
\widehat{z}}\cdot\sum_b\mathcal{U}_{ab}(M_p)\cdot\Theta_b(\widehat{Z},\widehat{z})
\end{eqnarray*}
where $\mathcal{U}(M_p)$ is a constant unitary $2^{n-1}\times 2^{n-1}$-matrix. So we
end up with the following first Fourier-coefficient
\begin{eqnarray*}
\sum_{p=1}^P (\det(\widehat{C}_p\widehat{Z}+\widehat{D}_p))^{-r}
g_1(T,\widehat{Z}_p,\widehat{z}_p)=
\sum_b \Psi(b,T)(\widehat{Z})\cdot\Theta_b(\widehat{Z},\widehat{z})
\end{eqnarray*}
with
\begin{eqnarray*}
\Psi(b,T)(\widehat{Z})=\sum_{a\neq(0,...,0),p=1,...,P}\mathcal{U}_{ab}(M_p)\cdot
(\det(\widehat{C}_p\widehat{Z}+\widehat{D}_p))^{-r+\frac{1}{2}}\cdot
\psi(a,T)(\widehat{Z}_p)
\end{eqnarray*}
for $T_{1n}=\frac{1}{2}$ and
\begin{eqnarray*}
\Psi(b,T)(\widehat{Z})=\sum_{p=1,...,P}\mathcal{U}_{(0),b}(M_p)\cdot
(\det(\widehat{C}_p\widehat{Z}+\widehat{D}_p))^{-r+\frac{1}{2}}\cdot
\psi((0,0,...,0),T)(\widehat{Z}_p)
\end{eqnarray*}
for $T_{1n}=0$. Now $\widehat{Z}_p=\widehat{X}_p+i\widehat{Y}_p$ and
\begin{eqnarray*}
\widehat{Y}_p=(\widehat{Z}\widehat{C}_p^t+\widehat{D}_p^t)^{-1}\cdot\widehat{Y}\cdot
(\widehat{C}_p\overline{\widehat{Z}}+\widehat{D}_p)^{-1}
\end{eqnarray*}
This implies because of $|\det(\widehat{C}_p\widehat{Z}+\widehat{D}_p)|=1$
\begin{eqnarray*}
\det(\widehat{Y}_p)=\det(\widehat{Y})
\end{eqnarray*}
In case $\widehat{Y}_p\neq\widehat{Y}$ for all $2\leq p\leq P$ we first choose
$\widehat{T}=\widehat{Y}^{-1}$ which is positive and symmetric. Then we have
$\sigma(\widehat{T}\widehat{Y})=n-1$ and
$\sigma(\widehat{T}\widehat{Y}_p)=
\sigma(\widehat{T}^{\frac{1}{2}}\widehat{Y}_p\widehat{T} ^{\frac{1}{2}})$. Let
\begin{eqnarray*}
\widetilde{Y}_p=\widehat{T}^{\frac{1}{2}}\widehat{Y}_p\widehat{T} ^{\frac{1}{2}}=
\widetilde{S}_p^t\widetilde{D}_p\widetilde{S}_p
\end{eqnarray*}
be the spectral decomposition with eigenvalues $\widetilde{D}_p=(\lambda_{p,kk}>0)$.
We have
\begin{eqnarray*}
\prod_{k=1}^{n-1}\lambda_{p,kk}=\det(\widetilde{D}_p)=\det(\widetilde{Y}_p)=
\det(\widehat{T})\cdot\det(\widehat{Y}_p)=\frac{det(\widehat{Y}_p)}{det(\widehat{Y})}=1
\end{eqnarray*}
Now not all $\lambda_{p,kk}$ can be equal because otherwise we would have
$\lambda_{p,kk}=1$ and so $\widetilde{Y}_p=\widehat{E}$ and so
$\widehat{Y}_p=\widehat{T}^{-1}=\widehat{Y}$ which was excluded. So by
arithmetic-geometric inequality we have
\begin{eqnarray*}
1=\left(\prod_{k=1}^{n-1}\lambda_{p,kk}\right)^{\frac{1}{n-1}}<
\frac{1}{n-1}\sum_{k=1}^{n-1}\lambda_{p,kk}=
\frac{1}{n-1}\sigma(\widetilde{Y}_p)=\frac{1}{n-1}\sigma(\widehat{T}\widehat{Y}_p)
\end{eqnarray*}
This shows that all positive symmetric $(n-1)\times(n-1)$-matrices with the same
determinant as $\widehat{Y}$ have coupled with 
$\widehat{T}=\widehat{Y}^{-1}$ a larger trace. Now choose a rational, symmetric matrix
$\widehat{T}$ which is sufficiently close to $\widehat{Y}^{-1}$ and scale by a large
number to make $\widehat{T}$ integer. Then all symplectic transformations different
from identity will make up a larger trace. The same argument as before now yields
$2^{n-1}$ matrices $T_1,...,T_{2^{n-1}}$ based on $\widehat{T}$ whose coefficient matrix
$(\Psi(b,T_j)(\widehat{Z}))_{b,j}$ has maximal rank $=2^{n-1}$. The case where some
$\widehat{Y}_p=\widehat{Y}$ but then $\widehat{X}_p\neq\widehat{X}$ mod $\mathbb{Z}$
can be handled as before. This shows the Theorem.\\
\\
For the proof of the Theorem we made use of the following fact
\begin{proposition}
For $M\in\Gamma_n$ with $C_{in}\neq 0$ for at least one $1\leq i\leq n$ we have
\begin{eqnarray*}
|\det(CZ+D)|\longrightarrow\infty\mbox{ for }y_{nn}\longrightarrow\infty
\end{eqnarray*} 
\end{proposition}
\textbf{Proof:} The case where $\det(C)\neq 0$ is easy and omitted. So let $\det(C)=0$
and $(C_{1n},C_{2n},...,C_{nn})\neq(0,0,...,0)$. Because of $\det(C)=0$ the columns of $C$
are $\mathbb{Q}$-linear dependent and so there are $u_1,...,u_n\in\mathbb{Z}$ with
no common divisor and
\begin{eqnarray*}
\sum_{m=1}^n u_m\cdot(C_{m1},...,C_{mn})=(0,0,...,0)
\end{eqnarray*}
Let $U\in GL(n,\mathbb{Z})$ with $U_{nm}=u_m$. Then the last row of $UC$ is
identical 0 and $|\det(CZ+D)|=|\det(UCZ+UD)|$. In the same way one can reduce $C$ by
left multiplication with unimodular matrices and not changing $|\det(CZ+D)|$ such that
\begin{eqnarray*}
C=\left(
\begin{array}{cc}
(C)_{jj},(C)_{j,n-j}\\
(0)_{n-j,j},(0)_{n-j,n-j}
\end{array}
\right)\mbox{ with }rank((C)_{j,n})=j
\end{eqnarray*}
and $(C_{1n},C_{2n},...,C_{jn})\neq(0,0,...,0)$. Let $V\in GL(n,\mathbb{Z})$ with
\begin{eqnarray*}
\widetilde{C}=C\cdot V=\left(
\begin{array}{cc}
(\widetilde{C})_{jj},0\\
0,0
\end{array}
\right)\mbox{ with }\det((\widetilde{C})_{j,j})\neq 0
\end{eqnarray*}
and let $\widetilde{Z}=WZW^t$ and $\widetilde{D}=DW^t$ with $W=V^{-1}$. Then
$|\det(CZ+D)|=|\det(\widetilde{C}\widetilde{Z}+\widetilde{D})|=
|\det((\widetilde{C})_{jj}\cdot(\widetilde{Z})_{jj}+(\widetilde{D})_{jj})|$. From
$C=\widetilde{C}W$ we read off
\begin{eqnarray*}
\sum_{k=1}^j\widetilde{C}_{mk}\cdot W_{kn}=C_{mn}\mbox{ for }m=1,2,...,j
\end{eqnarray*}
Now $(\widetilde{C})_{jj}$ is invertible and $(C_{1n},C_{2n},...,C_{jn})\neq(0,0,...,0)$ and
so $(W_{1n},W_{2n},...,W_{jn})\neq(0,0,...,0)$. Now we have
\begin{eqnarray*}
|\det((\widetilde{C})_{jj}\cdot(\widetilde{Z})_{jj}+(\widetilde{D})_{jj})|\geq
|\det((\widetilde{Z})_{jj}+(\widetilde{C})_{jj}^{-1}\cdot(\widetilde{D})_{jj})|\geq
\det((\widetilde{Y})_{jj})
\end{eqnarray*}
and $(\widetilde{Y})_{jj}=(W)_{jn}\cdot Y\cdot (W)_{jn}^t$. Now $Y$ is Minkowski-reduced
$Y\geq\delta_n(Y_{kk})$
and all eigenvalues of $Y$ are positive bounded below. Because of
$rank((W)_{jn})=j$ all eigenvalues of $(\widetilde{Y})_{jj}$ are also positive bounded below
independent of $y_{nn}\longrightarrow\infty$ (but maybe dependent on $W$).
Let $v=(W_{1n},...,W_{jn})/\|(W_{1n},...,W_{jn})\|$, then
\begin{eqnarray*}
v^t(\widetilde{Y})_{jj}v=\frac{1}{\|(W_{1n},...,W_{jn})\|^2}\sum_{p,q=1}^j
\widetilde{Y}_{pq}W_{pn}W_{qn}\geq\\
\frac{\delta_n}{\|(W_{1n},...,W_{jn})\|^2}\sum_{p,q=1}^j\left(
\sum_{k=1}^n Y_{kk}W_{pk}W_{qk}\right)W_{pn}W_{qn}\geq\\
\frac{\delta_n}{\|(W_{1n},...,W_{jn})\|^2}\sum_{p,q=1}^j Y_{nn}W_{pn}^2 W_{qn}^2=
\delta_n Y_{nn}\|(W_{1n},...,W_{jn})\|^2
\end{eqnarray*}
This shows that one eigenvalue of $(\widetilde{Y})_{jj}$ turns to $\infty$ for
$y_{nn}\longrightarrow\infty$ and so
\begin{eqnarray*}
\det((\widetilde{Y})_{jj})\longrightarrow\infty\mbox{ for }y_{nn}\longrightarrow\infty
\end{eqnarray*}
This proves the claim.\\
\\
\textbf{Remark 1:} In case $\widehat{Z}\in\mathcal{F}_{n-1}$ has nontrivial stabilizer
the matrix of coefficient modular forms
$(\phi_{1,a}(F)(\widehat{Z}))_{a\in(\frac{1}{2}\mathbb{Z}/\mathbb{Z})^{n-1},
F\in[\Gamma,r]}$ cannot have full rank $2^{n-1}$. Let for example
$\widehat{Z}=\widehat{U}^t\widehat{Z}\widehat{U}$ with
$\widehat{U}$ unimodular and $b\equiv\widehat{U}a$ mod $\mathbb{Z}$ where
$a,b$ are thetacharacteristics and $a\neq b$. Let
\begin{eqnarray*}
U=\left(
\begin{array}{cc}
\widehat{U}&0\\
0&1
\end{array}
\right)\mbox{ and so }
U^t ZU=\left(
\begin{array}{cc}
\widehat{U}^t\widehat{Z}\widehat{U}&\widehat{U}^t\widehat{z}\\
\widehat{z}^t\widehat{U}&z_{nn}
\end{array}
\right)
\end{eqnarray*}
For every Siegel modular form $F$ we get $F(U^t ZU)=F(Z)$ and so for the first
Fourier-Jacobi-Form $\Phi_1(F)(\widehat{Z},\widehat{z})=
\Phi_1(F)(\widehat{U}^t\widehat{Z}\widehat{U},\widehat{U}^t\widehat{z})$. This implies
\begin{eqnarray*}
\sum_{a\in(\frac{1}{2}\mathbb{Z}/\mathbb{Z})^{n-1}}
\phi_{1,a}(F)(\widehat{Z})\cdot\Theta_{1,a}(\widehat{Z},\widehat{z})=
\sum_{a\in(\frac{1}{2}\mathbb{Z}/\mathbb{Z})^{n-1}}
\phi_{1,a}(F)(\widehat{U}^t\widehat{Z}\widehat{U})\cdot
\Theta_{1,a}(\widehat{U}^t\widehat{Z}\widehat{U},\widehat{U}^t\widehat{z})
\end{eqnarray*}
Now with $b\equiv\widehat{U}a$ mod $\mathbb{Z}$ we get
\begin{eqnarray*}
\Theta_{1,a}(\widehat{U}^t\widehat{Z}\widehat{U},\widehat{U}^t\widehat{z})=
\sum_{k\in\mathbb{Z}^{n-1}}
e^{2\pi im(k+b)^t\widehat{Z}(k+b)+4\pi im(k+b)\widehat{z}}=
\Theta_{1,b}(\widehat{Z},\widehat{z})
\end{eqnarray*}
Using $\widehat{Z}=\widehat{U}^t\widehat{Z}\widehat{U}$ we have
$\phi_{1,a}(F)(\widehat{Z})=\phi_{1,b}(F)(\widehat{Z})$ for all Siegel modular forms $F$
and the special point $\widehat{Z}$. So the rank must be reduced in this case.\\
\\
\textbf{Remark 2:} Choosing $\widehat{U}=-\widehat{E}$ shows that 
$\Phi_m(F)(\widehat{Z},-\widehat{z})=\Phi_m(F)(\widehat{Z},\widehat{z})$, so only
even decompositions of $\Phi_m(F)(\widehat{Z},\widehat{z})$ in Thetaseries can occur.
For $m=1$ this is automatically fulfilled as $\Theta_{1,a}(\widehat{Z},\widehat{z})=
\Theta_{1,a}(\widehat{Z},-\widehat{z})$. The fiber over $\widehat{Z}$ to be embedded by
Siegel modular forms is therefore a Kummer-variety
$K(\widehat{Z})=(\mathbb{C}^{n-1}/\Lambda)/\{\pm 1\}$. The linear system
$(\Theta_{1,a}(\widehat{Z},\widehat{z}))_{a\in(\frac{1}{2}\mathbb{Z}/\mathbb{Z})^{n-1}}$
provides an embedding of the Kummer-variety \cite{7}. So by Theorem 1
there is a set of Siegel-cusp-forms combined with a set of Eisenstein-series which map
$(\mathcal{U}_{\varepsilon}(\widehat{Z})\times\mathbb{C}^{n-1}\times
(\mathbb{H}^1\cap\{y_{nn}>>1\}))/\Gamma_n$ biholomorphically onto
the total space of the Kummer-fibre-bundle. This is a local description of the Siegel
modular variety in a neighbourhood of a Satake boundary point.\\ 
\pagebreak\\
Now let $n=3$ and $F\in[\Gamma_3,r]$ be once again a Siegel cusp form. Then $F$
can be expanded
\begin{eqnarray*}
F(Z)=\Big(\sum_{a=(a_1,a_2),a_j\in\{0,\frac{1}{2}\}}\phi_a(\widehat{Z})\cdot
\Theta_a(\widehat{Z},\widehat{z})\Big)\cdot e^{2\pi iz_{33}}+
O\Big(e^{4\pi iz_{33}}\Big)
\end{eqnarray*}
with
\begin{eqnarray*}
\Theta_a(\widehat{Z},\widehat{z})=\sum_{(m_1,m_2)\in\mathbb{Z}^2}
e^{2\pi i((m_1+a_1)^2 z_{11}+2(m_1+a_1)(m_2+a_2)z_{12}+(m_2+a_2)^2 z_{22})+
4\pi i((m_1+a_1)z_{13}+(m_2+a_2)z_{23})}
\end{eqnarray*}
Let
\begin{eqnarray*}
S=\left(
\begin{array}{ccc}
S_{11}&S_{12}&0\\
S_{12}&S_{22}&0\\
0&0&0
\end{array}
\right)
\end{eqnarray*}
be a symmetric integer matrix. Then from $F(Z+S)=F(Z)$ and 
uniqueness of Fourier decomposition
we get
\begin{eqnarray*}
\sum_a\phi_a(\widehat{Z}+\widehat{S})\cdot
\Theta_a(\widehat{Z}+\widehat{S},\widehat{z})=
\sum_a\phi_a(\widehat{Z})\cdot\Theta_a(\widehat{Z},\widehat{z})
\end{eqnarray*}
Using $\Theta_a(\widehat{Z}+\widehat{S},\widehat{z})=
\Theta_a(\widehat{Z},\widehat{z})e^{2\pi ia^t\widehat{S}a}$ we get
\begin{eqnarray*}
\phi_a(\widehat{Z}+\widehat{S})=\phi_a(\widehat{Z})e^{2\pi ia^t\widehat{S}a}
\end{eqnarray*}
Let
\begin{eqnarray*}
U=\left(
\begin{array}{ccc}
1&1&0\\
0&1&0\\
0&0&1
\end{array}
\right)
\end{eqnarray*}
be a special unimodular matrix with action on a $Z$ by
\begin{eqnarray*}
U^t ZU=\left(
\begin{array}{ccc}
z_{11}&z_{11}+z_{12}&z_{13}\\
z_{11}+z_{12}&z_{11}+2z_{12}+z_{22}&z_{13}+z_{23}\\
z_{13}&z_{13}+z_{23}&z_{33}
\end{array}
\right)
\end{eqnarray*}
From $F(U^t ZU)=F(Z)$ we conclude
\begin{eqnarray*}
\sum_a\phi_a\left(
\begin{array}{cc}
z_{11}&z_{11}+z_{12}\\
z_{11}+z_{12}&z_{11}+2z_{12}+z_{22}
\end{array}
\right)
\Theta_a\left(\left(
\begin{array}{cc}
z_{11}&z_{11}+z_{12}\\
z_{11}+z_{12}&z_{11}+2z_{12}+z_{22}
\end{array}
\right),z_{13},z_{13}+z_{23}\right)=
\end{eqnarray*}
\begin{eqnarray*}
=\sum_a\phi_a(\widehat{Z})\cdot\Theta_a(\widehat{Z},\widehat{z})
\end{eqnarray*}
An easy calculation shows
\begin{eqnarray*}
\Theta_{(a_1,a_2)}\left(\left(
\begin{array}{cc}
z_{11}&z_{11}+z_{12}\\
z_{11}+z_{12}&z_{11}+2z_{12}+z_{22}
\end{array}
\right),z_{13},z_{13}+z_{23}\right)=
\Theta_{(a_1+a_2,a_2)}(\widehat{Z},\widehat{z})
\end{eqnarray*}
and so by uniqueness of Theta-decomposition
\begin{eqnarray*}
\phi_{(a_1,0)}\left(
\begin{array}{cc}
z_{11}&z_{11}+z_{12}\\
z_{11}+z_{12}&z_{11}+2z_{12}+z_{22}
\end{array}
\right)=\phi_{(a_1,0)}(\widehat{Z})\\
\phi_{(a_1,\frac{1}{2})}\left(
\begin{array}{cc}
z_{11}&z_{11}+z_{12}\\
z_{11}+z_{12}&z_{11}+2z_{12}+z_{22}
\end{array}
\right)=\phi_{(a_1+\frac{1}{2},\frac{1}{2})}(\widehat{Z})
\end{eqnarray*}
In case of $a_2=0$ and setting $S_{11}=0$ and so $a^t Sa=0$ we get
\begin{eqnarray*}
\phi_{(a_1,0)}\left(\widehat{Z}+\left(
\begin{array}{cc}
0&S_{12}\\
S_{12}&S_{22}
\end{array}
\right)\right)=\phi_{(a_1,0)}(\widehat{Z})
\end{eqnarray*}
So $\phi_{(a_1,0)}(\widehat{Z})$ is periodic in $z_{12},z_{22}$ and so has a Fourier series
\begin{eqnarray*}
\phi_{(a_1,0)}(\widehat{Z})=\sum_{k=1}^{\infty}\psi_{(a_1,0),k}(z_{11},z_{12})
e^{2\pi ikz_{22}}
\end{eqnarray*}
where the coefficient functions fulfill
\begin{eqnarray*}
\psi_{(a_1,0),k}(z_{11},z_{12}+1)=\psi_{(a_1,0),k}(z_{11},z_{12})\\
\psi_{(a_1,0),k}(z_{11},z_{12}+z_{11})e^{2\pi ik(z_{11}+2z_{12})}=
\psi_{(a_1,0),k}(z_{11},z_{12})
\end{eqnarray*}
So the function $\psi_{(a_1,0),1}(z_{11},z_{12})$ can be decomposed further by 1-dimensional
theta-functions
\begin{eqnarray*}
\psi_{(a_1,0),1}(z_{11},z_{12})=\psi_{(a_1,0),1,1}(z_{11})\cdot\Theta_{0}(z_{11},z_{12})+
\psi_{(a_1,0),1,2}(z_{11})\cdot\Theta_{\frac{1}{2}}(z_{11},z_{12})
\end{eqnarray*}
Here we set
\begin{eqnarray*}
\Theta_a(z_{11},z_{12})=\sum_{m\in\mathbb{Z}}
e^{2\pi i(m+a)^2 z_{11}+4\pi i(m+a)z_{12}}
\end{eqnarray*}
So in case of $a_2=0$ we have a decomposition
\begin{eqnarray*}
\phi_{(a_1,0)}(\widehat{Z})=\Big(\psi_{(a_1,0),1,1}(z_{11})\cdot\Theta_{0}(z_{11},z_{12})+
\psi_{(a_1,0),1,2}(z_{11})\cdot\Theta_{\frac{1}{2}}(z_{11},z_{12})\Big)\cdot
e^{2\pi iz_{22}}+e^{4\pi iz_{22}}\cdot...
\end{eqnarray*}
Similarly the corresponding Thetafunctions have an expansion of the form
\begin{eqnarray*}
\Theta_{(a_1,0)}(\widehat{Z},\widehat{z})=\Theta_{a_1}(z_{11},z_{13})+e^{2\pi iz_{22}}\cdot...
\end{eqnarray*}
which ends up in an expansion
\begin{eqnarray*}
\phi_{(a_1,0)}(\widehat{Z})\cdot\Theta_{(a_1,0)}(\widehat{Z},\widehat{z})=
\end{eqnarray*}
\begin{eqnarray*}
\Big(\psi_{(a_1,0),1,1}(z_{11})\cdot\Theta_{0}(z_{11},z_{12})+
\psi_{(a_1,0),1,2}(z_{11})\cdot\Theta_{\frac{1}{2}}(z_{11},z_{12})\Big)\cdot
\Theta_{a_1}(z_{11},z_{13})\cdot
e^{2\pi iz_{22}}+e^{4\pi iz_{22}}\cdot...
\end{eqnarray*}
In case of $a_2=\frac{1}{2}$ and $S_{11}=0$ we get
\begin{eqnarray*}
\phi_{(0,\frac{1}{2})}\left(\widehat{Z}+\left(
\begin{array}{cc}
0&S_{12}\\
S_{12}&S_{22}
\end{array}
\right)\right)=\phi_{(0,\frac{1}{2})}(\widehat{Z})\cdot e^{-\frac{1}{2}\pi iS_{22}}\\
\phi_{(\frac{1}{2},\frac{1}{2})}\left(\widehat{Z}+\left(
\begin{array}{cc}
0&S_{12}\\
S_{12}&S_{22}
\end{array}
\right)\right)=\phi_{(\frac{1}{2},\frac{1}{2})}(\widehat{Z})\cdot e^{-\frac{1}{2}\pi iS_{22}-
\pi iS_{12}}
\end{eqnarray*}
So $\phi_{(a_1,\frac{1}{2})}(\widehat{Z})$ have a Fourier series
\begin{eqnarray*}
\phi_{(a_1,\frac{1}{2})}(\widehat{Z})=\sum_{k=1}^{\infty}\psi_{(a_1,\frac{1}{2}),k}(z_{11},z_{12})
e^{(2k-\frac{1}{2})\pi iz_{22}}
\end{eqnarray*}
where the coefficient functions fulfill
\begin{eqnarray*}
\psi_{(0,\frac{1}{2}),k}(z_{11},z_{12}+1)=\psi_{(0,\frac{1}{2}),k}(z_{11},z_{12})\\
\psi_{(\frac{1}{2},\frac{1}{2}),k}(z_{11},z_{12}+1)=-\psi_{(\frac{1}{2},\frac{1}{2}),k}(z_{11},z_{12})
\end{eqnarray*}
Furthermore from
\begin{eqnarray*}
\phi_{(a_1,\frac{1}{2})}\left(
\begin{array}{cc}
z_{11}&z_{11}+z_{12}\\
z_{11}+z_{12}&z_{11}+2z_{12}+z_{22}
\end{array}
\right)=\phi_{(a_1+\frac{1}{2},\frac{1}{2})}(\widehat{Z})
\end{eqnarray*}
we get
\begin{eqnarray*}
\psi_{(a_1,\frac{1}{2}),k}(z_{11},z_{12}+z_{11})e^{(2k-\frac{1}{2})\pi i(z_{11}+2z_{12})}=
\psi_{(a_1+\frac{1}{2},\frac{1}{2}),k}(z_{11},z_{12})
\end{eqnarray*}
As before we are especially interested in $k=1$. Iterating the above relation twice we get
\begin{eqnarray*}
\psi_{(a_1,\frac{1}{2}),1}(z_{11},z_{12}+2z_{11})=
\psi_{(a_1,\frac{1}{2}),1}(z_{11},z_{12})\cdot e^{-6\pi iz_{11}-6\pi iz_{12}}
\end{eqnarray*}
From this we can derive a representation
\begin{eqnarray*}
\psi_{(0,\frac{1}{2}),1}(z_{11},z_{12})=\psi_{(0,\frac{1}{2}),1,1}(z_{11})\cdot
\sum_{q\in\mathbb{Z}}e^{\frac{2}{3}\pi i(3q)^2 z_{11}+2\pi i(3q)z_{12}}+\\
\psi_{(0,\frac{1}{2}),1,2}(z_{11})\cdot
\sum_{q\in\mathbb{Z}}e^{\frac{2}{3}\pi i(3q+1)^2 z_{11}+2\pi i(3q+1)z_{12}}+\\
\psi_{(0,\frac{1}{2}),1,3}(z_{11})\cdot
\sum_{q\in\mathbb{Z}}e^{\frac{2}{3}\pi i(3q+2)^2 z_{11}+2\pi i(3q+2)z_{12}}
\end{eqnarray*}
\begin{eqnarray*}
\psi_{(\frac{1}{2},\frac{1}{2}),1}(z_{11},z_{12})=\psi_{(\frac{1}{2},\frac{1}{2}),1,1}(z_{11})\cdot
\sum_{q\in\mathbb{Z}}e^{\frac{2}{3}\pi i((3q)^2+3q)z_{11}+\pi i(2\cdot 3q+1)z_{12}}+\\
\psi_{(\frac{1}{2},\frac{1}{2}),1,2}(z_{11})\cdot
\sum_{q\in\mathbb{Z}}e^{\frac{2}{3}\pi i((3q+1)^2+(3q+1))z_{11}+\pi i(2(3q+1)+1)z_{12}}+\\
\psi_{(\frac{1}{2},\frac{1}{2}),1,3}(z_{11})\cdot
\sum_{q\in\mathbb{Z}}e^{\frac{2}{3}\pi i((3q+2)^2+(3q+2))z_{11}+\pi i(2(3q+2)+1)z_{12}}
\end{eqnarray*}
or written in Theta-functions
\begin{eqnarray*}
\psi_{(0,\frac{1}{2}),1}(z_{11},z_{12})=\psi_{(0,\frac{1}{2}),1,1}(z_{11})\cdot
\Theta_{0}\Big(3z_{11},\frac{3}{2}z_{12}\Big)+\\
\psi_{(0,\frac{1}{2}),1,2}(z_{11})\cdot
\Theta_{\frac{1}{3}}\Big(3z_{11},\frac{3}{2}z_{12}\Big)+\\
\psi_{(0,\frac{1}{2}),1,3}(z_{11})\cdot
\Theta_{\frac{2}{3}}\Big(3z_{11},\frac{3}{2}z_{12}\Big)
\end{eqnarray*}
\begin{eqnarray*}
\psi_{(\frac{1}{2},\frac{1}{2}),1}(z_{11},z_{12})=\psi_{(\frac{1}{2},\frac{1}{2}),1,1}(z_{11})\cdot
\Theta_{\frac{1}{6}}\Big(3z_{11},\frac{3}{2}z_{12}\Big)\cdot e^{-\frac{1}{6}\pi iz_{11}}+\\
\psi_{(\frac{1}{2},\frac{1}{2}),1,2}(z_{11})\cdot
\Theta_{\frac{3}{6}}\Big(3z_{11},\frac{3}{2}z_{12}\Big)\cdot e^{-\frac{1}{6}\pi iz_{11}}+\\
\psi_{(\frac{1}{2},\frac{1}{2}),1,3}(z_{11})\cdot
\Theta_{\frac{5}{6}}\Big(3z_{11},\frac{3}{2}z_{12}\Big)\cdot e^{-\frac{1}{6}\pi iz_{11}}
\end{eqnarray*}
There are six coefficient-functions in $z_{11}$ occuring but due to the relation
\begin{eqnarray*}
\psi_{(a_1,\frac{1}{2}),1}(z_{11},z_{12}+z_{11})e^{\frac{3}{2}\pi i(z_{11}+2z_{12})}=
\psi_{(a_1+\frac{1}{2},\frac{1}{2}),1}(z_{11},z_{12})
\end{eqnarray*}
only three can be chosen independently. In fact we have the relations
\begin{eqnarray*}
\psi_{(\frac{1}{2},\frac{1}{2}),1,1}(z_{11})=\psi_{(0,\frac{1}{2}),1,3}(z_{11})\cdot
e^{\frac{1}{6}\pi iz_{11}}\\
\psi_{(\frac{1}{2},\frac{1}{2}),1,2}(z_{11})=\psi_{(0,\frac{1}{2}),1,1}(z_{11})\cdot
e^{\frac{1}{6}\pi iz_{11}}\\
\psi_{(\frac{1}{2},\frac{1}{2}),1,3}(z_{11})=\psi_{(0,\frac{1}{2}),1,2}(z_{11})\cdot
e^{\frac{1}{6}\pi iz_{11}}
\end{eqnarray*}
There are further restrictions coming from the unimodular transformation
\begin{eqnarray*}
U=\left(
\begin{array}{ccc}
1&0&0\\
0&-1&0\\
0&0&1
\end{array}
\right)
\end{eqnarray*}
which implies together with $F(U^t ZU)=F(Z)$ the further relations
\begin{eqnarray*}
\psi_{(a_1,a_2),1}(z_{11},-z_{12})=\psi_{(a_1,a_2),1}(z_{11},z_{12})
\end{eqnarray*}
Here we used that
\begin{eqnarray*}
\Theta_{(a_1,a_2)}\left(\left(
\begin{array}{cc}
z_{11}&-z_{12}\\
-z_{12}&z_{22}
\end{array}
\right),z_{13},-z_{23}\right)=
\Theta_{(a_1,a_2)}\left(\left(
\begin{array}{cc}
z_{11}&z_{12}\\
z_{12}&z_{22}
\end{array}
\right),z_{13},z_{23}\right)
\end{eqnarray*}
We apply the decomposition of $\psi_{(a_1,a_2),1}(z_{11},z_{12})$ into 1-dim.
Theta-functions. In case of $a_2=0$ this does not yield any further restrictions but
in case of $a_2=\frac{1}{2}$ we have
\begin{eqnarray*}
\Theta_0(3z_{11},-\frac{3}{2}z_{12})=\Theta_0(3z_{11},\frac{3}{2}z_{12})\\
\Theta_{\frac{3}{6}}(3z_{11},-\frac{3}{2}z_{12})=\Theta_{\frac{3}{6}}(3z_{11},\frac{3}{2}z_{12})\\
\Theta_{\frac{1}{3}}(3z_{11},-\frac{3}{2}z_{12})=\Theta_{\frac{2}{3}}(3z_{11},-\frac{3}{2}z_{12})\\
\Theta_{\frac{1}{6}}(3z_{11},-\frac{3}{2}z_{12})=\Theta_{\frac{5}{6}}(3z_{11},-\frac{3}{2}z_{12})
\end{eqnarray*}
From uniqueness of Theta-decomposition this implies
\begin{eqnarray*}
\psi_{(0,\frac{1}{2}),1,2}(z_{11})=\psi_{(0,\frac{1}{2}),1,3}(z_{11})
\end{eqnarray*}
So in fact there are only two independent coefficient functions 
$\psi_{(0,\frac{1}{2}),1,1}(z_{11}),\psi_{(0,\frac{1}{2}),1,2}(z_{11})$ in case of $a_2=\frac{1}{2}$.\\
We develop the two-dimensional Theta-functions as well and obtain
\begin{eqnarray*}
\Theta_{(a_1,\frac{1}{2})}(\widehat{Z},\widehat{z})=e^{\frac{1}{2}\pi iz_{22}}\cdot\Big(
\sum_{m_1}e^{2\pi i((m_1+a_1)^2 z_{11}+(m_1+a_1)z_{12})+4\pi i((m_1+a_1)z_{13}+
\frac{1}{2}z_{23})}+\\
\sum_{m_1}e^{2\pi i((m_1+a_1)^2 z_{11}-(m_1+a_1)z_{12})+4\pi i((m_1+a_1)z_{13}-
\frac{1}{2}z_{23})}\Big)+e^{\frac{9}{2}\pi iz_{22}}\cdot...=\\
e^{\frac{1}{2}\pi iz_{22}}\cdot\Big(\Theta_{a_1}(z_{11},z_{13}+\frac{1}{2}z_{12})
\cdot e^{2\pi iz_{23}}+\Theta_{a_1}(z_{11},z_{13}-\frac{1}{2}z_{12})
\cdot e^{-2\pi iz_{23}}\Big)+e^{\frac{9}{2}\pi iz_{22}}\cdot...
\end{eqnarray*} 
Putting all expansions together we obtain an expansion for the Fourier-Jacobi cusp form
of index 1
\begin{eqnarray*}
\sum_{a=(a_1,a_2),a_j\in\{0,\frac{1}{2}\}}\phi_a(\widehat{Z})\cdot
\Theta_a(\widehat{Z},\widehat{z})=e^{2\pi iz_{22}}\cdot\Delta(z_{11},z_{12},z_{13},z_{23})+
e^{4\pi iz_{22}}\cdot...
\end{eqnarray*}
where
\begin{eqnarray*}
\Delta(z_{11},z_{12},z_{13},z_{23})=
\end{eqnarray*}
\begin{eqnarray*}
=\Big(\psi_{(0,0),1,1}(z_{11})\cdot\Theta_{0}(z_{11},z_{12})+
\psi_{(0,0),1,2}(z_{11})\cdot\Theta_{\frac{1}{2}}(z_{11},z_{12})\Big)\cdot
\Theta_{0}(z_{11},z_{13})
+\\
\Big(\psi_{(\frac{1}{2},0),1,1}(z_{11})\cdot\Theta_{0}(z_{11},z_{12})+
\psi_{(\frac{1}{2},0),1,2}(z_{11})\cdot\Theta_{\frac{1}{2}}(z_{11},z_{12})\Big)\cdot
\Theta_{\frac{1}{2}}(z_{11},z_{13})+\\
\Big(\psi_{(0,\frac{1}{2}),1,1}(z_{11})\cdot
\Theta_{0}(3z_{11},\frac{3}{2}z_{12})+
\psi_{(0,\frac{1}{2}),1,2}(z_{11})\cdot\\
\Big(
\Theta_{\frac{1}{3}}(3z_{11},\frac{3}{2}z_{12})+
\Theta_{\frac{2}{3}}(3z_{11},\frac{3}{2}z_{12})\Big)\Big)\cdot\\
\Big(\Theta_{0}(z_{11},z_{13}+\frac{1}{2}z_{12})
\cdot e^{2\pi iz_{23}}+\Theta_{0}(z_{11},z_{13}-\frac{1}{2}z_{12})
\cdot e^{-2\pi iz_{23}}\Big)+\\
\Big(\psi_{(0,\frac{1}{2}),1,1}(z_{11})\cdot
\Theta_{\frac{3}{6}}(3z_{11},\frac{3}{2}z_{12})+
\psi_{(0,\frac{1}{2}),1,2}(z_{11})\cdot\\
\Big(\Theta_{\frac{1}{6}}(3z_{11},\frac{3}{2}z_{12})+
\Theta_{\frac{5}{6}}(3z_{11},\frac{3}{2}z_{12})\Big)\Big)\cdot\\
\Big(\Theta_{\frac{1}{2}}(z_{11},z_{13}+\frac{1}{2}z_{12})
\cdot e^{2\pi iz_{23}}+\Theta_{\frac{1}{2}}(z_{11},z_{13}-\frac{1}{2}z_{12})
\cdot e^{-2\pi iz_{23}}\Big)
\end{eqnarray*}
Now the interesting point is that this expansion also holds true more generally for all
Jacobi cusp forms of scalar index 1 on $\mathbb{H}_2\times\mathbb{C}^2$.
But in case of a Fourier-Jacobi cusp form so stemming from a Siegel cusp form $F(Z)$ on
$\mathbb{H}_3$ we have an additional symmetry because of the unimodular identity
\begin{eqnarray*}
F\left(
\begin{array}{ccc}
z_{11}&z_{12}&z_{13}\\
z_{12}&z_{22}&z_{23}\\
z_{13}&z_{23}&z_{33}
\end{array}
\right)=
F\left(
\begin{array}{ccc}
z_{11}&z_{13}&z_{12}\\
z_{13}&z_{33}&z_{23}\\
z_{12}&z_{23}&z_{22}
\end{array}
\right)
\end{eqnarray*}
This additional symmetry implies
\begin{eqnarray*}
\Delta(z_{11},z_{12},z_{13},z_{23})=\Delta(z_{11},z_{13},z_{12},z_{23})
\end{eqnarray*}
As the torus variable $z_{23}$ may vary free this 
translates into the following extra conditions on the coefficient
function $\psi_{...}(z_{11})$:
\begin{eqnarray*}
\big(\psi_{(0,0),1,1}(z_{11})\cdot\Theta_{0}(z_{11},z_{12})+
\psi_{(0,0),1,2}(z_{11})\cdot\Theta_{\frac{1}{2}}(z_{11},z_{12})\big)\cdot
\Theta_{0}(z_{11},z_{13})
+\\
\big(\psi_{(\frac{1}{2},0),1,1}(z_{11})\cdot\Theta_{0}(z_{11},z_{12})+
\psi_{(\frac{1}{2},0),1,2}(z_{11})\cdot\Theta_{\frac{1}{2}}(z_{11},z_{12})\big)\cdot
\Theta_{\frac{1}{2}}(z_{11},z_{13})=
\end{eqnarray*}
\begin{eqnarray*}
=\big(\psi_{(0,0),1,1}(z_{11})\cdot\Theta_{0}(z_{11},z_{13})+
\psi_{(0,0),1,2}(z_{11})\cdot\Theta_{\frac{1}{2}}(z_{11},z_{13})\big)\cdot
\Theta_{0}(z_{11},z_{12})
+\\
\big(\psi_{(\frac{1}{2},0),1,1}(z_{11})\cdot\Theta_{0}(z_{11},z_{13})+
\psi_{(\frac{1}{2},0),1,2}(z_{11})\cdot\Theta_{\frac{1}{2}}(z_{11},z_{13})\big)\cdot
\Theta_{\frac{1}{2}}(z_{11},z_{12})
\end{eqnarray*}
and for the $a_2=\frac{1}{2}$ part
\begin{eqnarray*}
\Big(\psi_{(0,\frac{1}{2}),1,1}(z_{11})\cdot
\Theta_{0}(3z_{11},\frac{3}{2}z_{12})+
\psi_{(0,\frac{1}{2}),1,2}(z_{11})\cdot
\Theta_{\frac{1}{3}}(3z_{11},\frac{3}{2}z_{12})+\\
\psi_{(0,\frac{1}{2}),1,2}(z_{11})\cdot
\Theta_{\frac{2}{3}}(3z_{11},\frac{3}{2}z_{12})\Big)\cdot
\Theta_{0}(z_{11},z_{13}\pm\frac{1}{2}z_{12})+\\
\Big(\psi_{(0,\frac{1}{2}),1,1}(z_{11})\cdot
\Theta_{\frac{3}{6}}(3z_{11},\frac{3}{2}z_{12})+
\psi_{(0,\frac{1}{2}),1,2}(z_{11})\cdot
\Theta_{\frac{1}{6}}(3z_{11},\frac{3}{2}z_{12})+\\
\psi_{(0,\frac{1}{2}),1,2}(z_{11})\cdot
\Theta_{\frac{5}{6}}(3z_{11},\frac{3}{2}z_{12})\Big)\cdot
\Theta_{\frac{1}{2}}(z_{11},z_{13}\pm\frac{1}{2}z_{12})=
\end{eqnarray*}
\begin{eqnarray*}
=\Big(\psi_{(0,\frac{1}{2}),1,1}(z_{11})\cdot
\Theta_{0}(3z_{11},\frac{3}{2}z_{13})+
\psi_{(0,\frac{1}{2}),1,2}(z_{11})\cdot
\Theta_{\frac{1}{3}}(3z_{11},\frac{3}{2}z_{13})+\\
\psi_{(0,\frac{1}{2}),1,2}(z_{11})\cdot
\Theta_{\frac{2}{3}}(3z_{11},\frac{3}{2}z_{13})\Big)\cdot
\Theta_{0}(z_{11},z_{12}\pm\frac{1}{2}z_{13})+\\
\Big(\psi_{(0,\frac{1}{2}),1,1}(z_{11})\cdot
\Theta_{\frac{3}{6}}(3z_{11},\frac{3}{2}z_{13})+
\psi_{(0,\frac{1}{2}),1,3}(z_{11})\cdot
\Theta_{\frac{1}{6}}(3z_{11},\frac{3}{2}z_{13})+\\
\psi_{(0,\frac{1}{2}),1,2}(z_{11})\cdot
\Theta_{\frac{5}{6}}(3z_{11},\frac{3}{2}z_{13})\Big)\cdot
\Theta_{\frac{1}{2}}(z_{11},z_{12}\pm\frac{1}{2}z_{13})
\end{eqnarray*}
The first equation yields
\begin{eqnarray*}
(\psi_{(0,0),1,2}(z_{11})-\psi_{(\frac{1}{2},0),1,1}(z_{11}))\cdot\Theta_{\frac{1}{2}}(z_{11},z_{12})
\cdot\Theta_0(z_{11},z_{13})=\\
(\psi_{(0,0),1,2}(z_{11})-\psi_{(\frac{1}{2},0),1,1}(z_{11}))\cdot\Theta_{\frac{1}{2}}(z_{11},z_{13})
\cdot\Theta_0(z_{11},z_{12})
\end{eqnarray*}
As the torus parameter $z_{12},z_{13}$ may vary free we conclude from
linear independence of the theta-functions
\begin{eqnarray*}
\psi_{(0,0),1,2}(z_{11})=\psi_{(\frac{1}{2},0),1,1}(z_{11})
\end{eqnarray*}
In the appendix we prove that
\begin{eqnarray*}
\Delta_{1,\pm}(z_{11},z_{12},z_{13})=\Theta_{0}(3z_{11},\frac{3}{2}z_{12})\cdot
\Theta_{0}(z_{11},z_{13}\pm\frac{1}{2}z_{12})+
\Theta_{\frac{3}{6}}(3z_{11},\frac{3}{2}z_{12})\cdot
\Theta_{\frac{1}{2}}(z_{11},z_{13}\pm\frac{1}{2}z_{12})
\end{eqnarray*}
and
\begin{eqnarray*}
\Delta_{2,\pm}(z_{11},z_{12},z_{13})=
\Big(\Theta_{\frac{1}{3}}(3z_{11},\frac{3}{2}z_{12})+
\Theta_{\frac{1}{3}}(3z_{11},\frac{3}{2}z_{12})\Big)\cdot
\Theta_{0}(z_{11},z_{13}\pm\frac{1}{2}z_{12})+\\
\Big(\Theta_{\frac{1}{6}}(3z_{11},\frac{3}{2}z_{12})+
\Theta_{\frac{5}{6}}(3z_{11},\frac{3}{2}z_{12})\Big)\cdot
\Theta_{\frac{1}{2}}(z_{11},z_{13}\pm\frac{1}{2}z_{12})
\end{eqnarray*}
are symmetric in $z_{12},z_{13}$. So symmetry of $z_{12},z_{13}$ does not yield any
further restrictions on $\psi_{(0,\frac{1}{2}),1,1}(z_{11}),\psi_{(0,\frac{1}{2}),1,2}(z_{11})$.\\
\\
Now we construct a Jacobi cusp form of index 1 and arbitrary high weight $r$ for
which the symmetry relation $\psi_{(0,0),1,2}(z_{11})=\psi_{(\frac{1}{2},0),1,1}(z_{11})$
is not fulfilled. Let $J_{2,1}\subset\Gamma_3$ be the Jacobi-subgroup. Let
\begin{eqnarray*}
T=\left(
\begin{array}{ccc}
1&\frac{1}{2}&0\\
\frac{1}{2}&1&0\\
0&0&1
\end{array}
\right)
\end{eqnarray*}
and
\begin{eqnarray*}
\widetilde{F}_r(T,Z)=\sum_{[M]\in J_{2,1}\backslash\{Transl.\}}\det(CZ+D)^{-r}\cdot
e^{2\pi i\sigma(TM(Z))}\\
=\widehat{F}_r(T,\widehat{Z},\widehat{z})\cdot e^{2\pi iz_{33}}
\end{eqnarray*}
with a Jacobi-cusp form $\widehat{F}_r(T,\widehat{Z},\widehat{z})$ of weight $r$
and scalar index 1. We fix a Minkowski reduced $Y>0$ with $1<<Y_{11}<Y_{22}<Y_{33}$ and
$Y_{11}$ so large that for all $C\neq(0)$ we have
$|\det(CZ+D)|>2$. Now we have
\begin{eqnarray*}
\phi_{(a_1,0)}(\widehat{F}_r)(\widehat{Z})=
\frac{1}{\Theta_{(a_1,0)}(2\widehat{Y},2\widehat{y})}\cdot e^{2\pi y_{33}}\cdot\\
\int_0^1\int_0^1\int_0^1\widetilde{F}_r(T,Z)\cdot\overline{
\Theta_{(a_1,0)}(\widehat{Z},\widehat{z})}\cdot e^{-2\pi ix_{33}}
dx_{13}dx_{23}dx_{33}
\end{eqnarray*}
and employing this expression we get
\begin{eqnarray*}
\psi_{(a_1,0),1,1}(\widehat{F}_r)(z_{11})=
\frac{1}{\Theta_{0}(2Y_{11},2Y_{12})}\cdot e^{2\pi y_{22}}\cdot\\
\int_0^1\int_0^1\phi_{(a_1,0)}(\widehat{F}_r)(\widehat{Z})
\cdot\overline{
\Theta_{0}(Z_{11},Z_{12})}\cdot e^{-2\pi ix_{22}}
dx_{12}dx_{22}
\end{eqnarray*}
and a similar expression for $\psi_{(a_1,0),1,2}(\widehat{F}_r)(z_{11})$. From this we see
that all terms of $\widetilde{F}_r(T,Z)$ with $C\neq(0)$ will have a vanishing contribution
to $\psi_{(a_1,0),1,1}(\widehat{F}_r)(z_{11})$, $\psi_{(a_1,0),1,2}(\widehat{F}_r)(z_{11})$
for $r\longrightarrow\infty$. So for disproving the equality
$\psi_{(0,0),1,2}(\widehat{F}_r)(z_{11})=\psi_{(\frac{1}{2},0),1,1}(\widehat{F}_r)(z_{11})$
it is sufficient to restrict to the unimodular subseries of $\widetilde{F}_r(T,Z)$ given by
\begin{eqnarray*}
\widetilde{G}(T,Z)=\sum_U e^{2\pi i\sigma(UTU^t Z)}
\end{eqnarray*}
where
\begin{eqnarray*}
U=\left(
\begin{array}{ccc}
U_{11}&U_{12}&m_1\\
U_{21}&U_{22}&m_2\\
0&0&1
\end{array}
\right)\mbox{ and }\widehat{U}=\left(
\begin{array}{cc}
U_{11}&U_{12}\\
U_{21}&U_{22}
\end{array}
\right)\in GL(2,\mathbb{Z})
\end{eqnarray*}
and $m=(m_1,m_2)\in\mathbb{Z}^2$ vary free. So from
\begin{eqnarray*}
\sigma(UTU^t Z)=\sigma(\widehat{U}\widehat{T}\widehat{U}^t\widehat{Z})+
m^t\widehat{Z}m+2m^t\widehat{z}+z_{33}
\end{eqnarray*}
we immediately obtain
\begin{eqnarray*}
\widetilde{G}(T,Z)=\left(\sum_{\widehat{U}\in GL(2,\mathbb{Z})}
e^{2\pi i\sigma(\widehat{U}\widehat{T}\widehat{U}^t\widehat{Z})}\right)\cdot
e^{2\pi iz_{33}}\cdot\Theta_{(0,0)}(\widehat{Z},\widehat{z})
\end{eqnarray*}
Furthermore a direct calculation shows
\begin{eqnarray*}
\sum_{\widehat{U}\in GL(2,\mathbb{Z})}
e^{2\pi i\sigma(\widehat{U}\widehat{T}\widehat{U}^t\widehat{Z})}=
12\cdot e^{-\frac{1}{2}\pi iz_{11}+2\pi iz_{22}}\cdot\Theta_{\frac{1}{2}}(z_{11},z_{12})+
e^{4\pi iz_{22}}\cdot...
\end{eqnarray*}
and so
\begin{eqnarray*}
\psi_{(0,0),1,2}(\widetilde{G})(z_{11})=12\cdot e^{-\frac{1}{2}\pi iz_{11}}\mbox{ and }
\psi_{(\frac{1}{2},0),1,1}(\widetilde{G})(z_{11})=0
\end{eqnarray*}
This proves
$\psi_{(0,0),1,2}(\widehat{F}_r)(z_{11})\neq\psi_{(\frac{1}{2},0),1,1}(\widehat{F}_r)(z_{11})$
for all large $r$.
\begin{theorem}
For all large weights $r$ there are more Jacobi cusp forms of scalar index 1 than Fourier-Jacobi cusp forms of index 1.
\end{theorem}
\small
Authors address: Bert Koehler, Debeka-Hauptverwaltung Ferdinand-Sauerbruch-Str 18,
56058 Koblenz, Email: Bert.Koehler@debeka.de\\
\\

\begin{center}
\large
\textbf{Appendix}
\end{center}
First we notice that from symmetry of $\Delta_{1,+}$
we immediately also obtain the corresponding symmetry for $\Delta_{1,-}$ as
\begin{eqnarray*}
\Delta_{1,+}(z_{11},-z_{12},z_{13})=\Theta_{0}(3z_{11},-\frac{3}{2}z_{12})\cdot
\Theta_{0}(z_{11},z_{13}-\frac{1}{2}z_{12})+\\
\Theta_{\frac{3}{6}}(3z_{11},-\frac{3}{2}z_{12})\cdot
\Theta_{\frac{1}{2}}(z_{11},z_{13}-\frac{1}{2}z_{12})=\\
\Theta_{0}(3z_{11},\frac{3}{2}z_{12})\cdot
\Theta_{0}(z_{11},z_{13}-\frac{1}{2}z_{12})+\\
\Theta_{\frac{3}{6}}(3z_{11},\frac{3}{2}z_{12})\cdot
\Theta_{\frac{1}{2}}(z_{11},z_{13}-\frac{1}{2}z_{12})=\Delta_{1,-}(z_{11},z_{12},z_{13})
\end{eqnarray*}
and so
\begin{eqnarray*}
\Delta_{1,-}(z_{11},z_{13},z_{12})=\Delta_{1,+}(z_{11},-z_{13},z_{12})=
\Delta_{1,+}(z_{11},z_{12},-z_{13})=\\
\Delta_{1,+}(z_{11},-z_{12},z_{13})=
\Delta_{1,-}(z_{11},z_{12},z_{13})
\end{eqnarray*}
A similar calculation holds for $\Delta_{2,-}$. For $b=0,\frac{1}{3},\frac{2}{3}$ we set
by abuse of notation
\begin{eqnarray*}
\Delta_b(z_{11},z_{12},z_{13})=\Theta_b(3z_{11},\frac{3}{2}z_{12})\cdot
\Theta_0(z_{11},z_{13}+\frac{1}{2}z_{12})+
\Theta_{b+\frac{1}{2}}(3z_{11},\frac{3}{2}z_{12})\cdot
\Theta_{\frac{1}{2}}(z_{11},z_{13}+\frac{1}{2}z_{12})
\end{eqnarray*}
We claim
\begin{eqnarray*}
(z_{12},z_{13})\longmapsto\Delta_b(z_{11},z_{12},z_{13})
\end{eqnarray*}
is a symmetric function.\\
\\
We have
\begin{eqnarray*}
\Theta_b(3z_{11},\frac{3}{2}z_{12})=\sum_{q\in\mathbb{Z}} e^{2\pi i(q+b)^2\cdot 3z_{11}+
2\pi i(q+b)\cdot 3z_{12}}=\\
\sum_{q\in\mathbb{Z}}\sum_{j=0}^1 e^{6\pi i(2q+j+b)^2 z_{11}+
6\pi i(2q+j+b)z_{12}}
\end{eqnarray*}
and
\begin{eqnarray*}
\Theta_0(z_{11},z_{13}+\frac{1}{2}z_{12})=
\sum_{m\in\mathbb{Z}}e^{2\pi im^2 z_{11}+4\pi im(z_{13}+\frac{1}{2}z_{12})}=\\
\sum_{m\in\mathbb{Z}}\sum_{k=0}^5 e^{2\pi i(6m+k)^2 z_{11}+2\pi i(6m+k)(2z_{13}+z_{12})}
\end{eqnarray*}
So multiplying we get
\begin{eqnarray*}
\Delta_b(z_{11},z_{12},z_{13})=
\sum_{q,m\in\mathbb{Z}}\sum_{j=0}^1\sum_{k=0}^5 
e^{2\pi i(3(2q+j+b)^2+(6m+k)^2)z_{11}+2\pi i(6q+3j+3b+6m+k)z_{12}+
2\pi i(12m+2k)z_{13}}\\
+\sum_{q,m\in\mathbb{Z}}\sum_{j=0}^1\sum_{k=0}^5 
e^{2\pi i(3(2q+j+b+\frac{1}{2})^2+(6m+k+\frac{1}{2})^2)z_{11}+
2\pi i(6q+3j+3b+\frac{3}{2}+6m+k+\frac{1}{2})z_{12}+
2\pi i(12m+2k+1)z_{13}}
\end{eqnarray*}
\begin{eqnarray*}
=\sum_{q,m\in\mathbb{Z}}\sum_{j=0}^1\sum_{k=0}^5 
e^{2\pi i(3(2q-2m+j+b)^2+(6m+k)^2)z_{11}+2\pi i(6q+3j+3b+k)z_{12}+
2\pi i(12m+2k)z_{13}}\\
+\sum_{q,m\in\mathbb{Z}}\sum_{j=0}^1\sum_{k=0}^5 
e^{2\pi i(3(2q-2m+j+b+\frac{1}{2})^2+(6m+k+\frac{1}{2})^2)z_{11}+
2\pi i(6q+3j+3b+2+k)z_{12}+
2\pi i(12m+2k+1)z_{13}}
\end{eqnarray*}
\begin{eqnarray*}
=\sum_{q,m\in\mathbb{Z}}\sum_{j,l=0}^1\sum_{k=0}^5 
e^{2\pi i(3(4q-2m+2l+j+b)^2+(6m+k)^2)z_{11}+2\pi i(12q+6l+3j+3b+k)z_{12}+
2\pi i(12m+2k)z_{13}}\\
+\sum_{q,m\in\mathbb{Z}}\sum_{j,l=0}^1\sum_{k=0}^5 
e^{2\pi i(3(4q-2m+2l+j+b+\frac{1}{2})^2+(6m+k+\frac{1}{2})^2)z_{11}+
2\pi i(12q+6l+3j+3b+k+2)z_{12}+2\pi i(12m+2k+1)z_{13}}
\end{eqnarray*}
For such an expression to be symmetric in $z_{12},z_{13}$ it is necessary that for every
choice of $j,l,k$ there is a corresponding tripel $\widetilde{j},\widetilde{l},\widetilde{k}$
with either
\begin{eqnarray*}
6l+3j+3b+k\equiv 2\widetilde{k}\mbox{ mod }12\\
6\widetilde{l}+3\widetilde{j}+3b+\widetilde{k}\equiv 2k\mbox{ mod }12
\end{eqnarray*}
or
\begin{eqnarray*}
6l+3j+3b+k\equiv 2\widetilde{k}+1\mbox{ mod }12\\
6\widetilde{l}+3\widetilde{j}+3b+\widetilde{k}+2\equiv 2k\mbox{ mod }12
\end{eqnarray*}
So let $j,l,k$ be fixed and assume first that 
$6l+3j+3b+k=2\widetilde{k}$ is even. Then we conclude
\begin{eqnarray*}
k+3b-2\widetilde{k}=-6l-3j\equiv 0\mbox{ mod }3\\
\Rightarrow 4k-6b-2\widetilde{k}\equiv 0\mbox{ mod }3\\
\Rightarrow 2k-3b-\widetilde{k}\equiv 0\mbox{ mod }3\\
\Rightarrow\exists\widetilde{j},\widetilde{l}\mbox{ with }
6\widetilde{l}+3\widetilde{j}+3b+\widetilde{k}\equiv 2k\mbox{ mod }12
\end{eqnarray*}
In case 
$6l+3j+3b+k=2\widetilde{k}+1$ is odd we conclude
\begin{eqnarray*}
k+3b-2\widetilde{k}-1=-6l-3j\equiv 0\mbox{ mod }3\\
\Rightarrow 4k-6b-2\widetilde{k}-4\equiv 0\mbox{ mod }3\\
\Rightarrow 2k-3b-\widetilde{k}-2\equiv 0\mbox{ mod }3\\
\Rightarrow\exists\widetilde{j},\widetilde{l}\mbox{ with }
6\widetilde{l}+3\widetilde{j}+3b+\widetilde{k}+2\equiv 2k\mbox{ mod }12
\end{eqnarray*}
So for every pair of Fourier-coefficients
$(12q+6l+3j+3b+k,12m+2k)$ or $(12q+6l+3j+3b+2+k,12m+2k+1)$ 
of $z_{12},z_{13}$ there is a
corresponding pair $(12\widetilde{m}+2\widetilde{k},
12\widetilde{q}+6\widetilde{l}+3\widetilde{j}+3b+\widetilde{k})$ or
$(12\widetilde{m}+2\widetilde{k}+1,
12\widetilde{q}+6\widetilde{l}+3\widetilde{j}+3b+\widetilde{k}+2)$ 
with the same coefficients but reversed roles of $z_{12},z_{13}$.\\
In the next step we have to prove that for such corresponding pairs the
coefficients of $z_{11}$ are identical. So assume first we have
\begin{eqnarray*}
6l+3j+3b+k=2\widetilde{k}\mbox{ and}\\
6\widetilde{l}+3\widetilde{j}+3b+\widetilde{k}=2k
\end{eqnarray*}
We calculate the $z_{11}$-coefficient
\begin{eqnarray*}
3(4q-2m+2l+j+b)^2+(6m+k)^2=48q^2-48qm+48m^2+\\
12(4l+2j+2b)q+12(k-(2l+j+b))m+3(2l+j+b)^2+k^2
\end{eqnarray*}
So we have to show that
\begin{eqnarray*}
4l+2j+2b=\widetilde{k}-(2\widetilde{l}+\widetilde{j}+b)\\
4\widetilde{l}+2\widetilde{j}+2b=k-(2l+j+b)\\
3(2l+j+b)^2+k^2=3(2\widetilde{l}+\widetilde{j}+b)^2+\widetilde{k}^2
\end{eqnarray*}
But from our two relations above we infer
\begin{eqnarray*}
2l+j+b=\frac{2}{3}\widetilde{k}-\frac{1}{3}k\\
2\widetilde{l}+\widetilde{j}+b=\frac{2}{3}k-\frac{1}{3}\widetilde{k}
\end{eqnarray*}
Plugging these terms in we directly obtain equality of $z_{11}$-coefficients. The calculation
for the remaining cases is similar.


\begin{thebibliography}{1}
\bibitem[1]{1} Freitag, E., Siegelsche Modulfunktionen, Springer 1983\\
\bibitem[2]{2} Eichler, M. and Zagier, D., The Theory of Jacobi forms, Birkh\"auser 1985\\
\bibitem[3]{3} Dulinski, J., A decomposition theorem for Jacobi forms, Math. Ann. 303, No. 3, 473-498 (1995)\\
\bibitem[4]{4} Klingen, H., \"Uber Kernfunktionen f\"ur Jacobiformen und Siegelsche Modulformen, Math. Ann. 285, 405-416 (1989)\\
\bibitem[5]{5} Igusa, J.I., Theta Functions, Springer 1972\\
\bibitem[6]{6} Klingen, H., Zum Darstellungssatz f\"ur Siegelsche Modulformen, Math. Z. 102, 30-43 (1967)\\
\bibitem[7]{7} Birkenhake, Ch. and Lange, H., Complex Abelian varieties, Springer 1992
\end{thebibliography}
\end{document}